\documentclass[12pt]{amsart}

\usepackage{amssymb, amsmath, amsthm, mathtools,}
\usepackage{indentfirst}
\usepackage[all]{xy}
\usepackage{cite}
\usepackage{comment}

\theoremstyle{plain}
\newtheorem{prop} {Proposition} [section]
\newtheorem{thm}[prop] {Theorem} 
\newtheorem{dfn} [prop]{Definition}
\newtheorem{exam}[prop]{Example}
\newtheorem{cor}[prop]{Corollary}

\theoremstyle{definition}
\newtheorem{prob*}[prop]{Problem}
\newtheorem{rem}[prop]{Remark}

\title[SOME COMPACT-LIKE PROPERTIES]%
{Some compact-like properties in non-archimedean functional analysis}
\author[ISHIZUKA]%
{KOSUKE ISHIZUKA}
\address{Mathematical Institute, Graduate School of Science, Tohoku University, 6-3 Aramakiaza, Aoba, Sendai, Miyagi 980-8578, Japan.}
\date{}
\email{kosuke.ishizuka.r7@dc.tohoku.ac.jp}


\begin{document}

\begin{abstract}
First, we define some concepts similar to the local compactoidity or the c-compactness, and study relationships between these concepts and the original ones. As a result, we find a characterization of the local compactoidity when its coefficient field is spherically complete. Moreover, from the point of view of the minimum principle, we give a necessary and sufficient condition for the c-compactness under a suitable condition. Secondly, we try a new approach to a non-complete local compactoid, which gives us a different perspective than before. Thirdly, we study the non-archimedean Goldstine theorem and Eberlein-\v{S}mulian theorem. Consequently, if the coefficient field is spherically complete, we get results completely different from the classical ones. Finally, we give a new result about the closed range theorem by using epicompactness. 
\end{abstract}

\maketitle

\section{Introduction and preliminaries}

\subsection{Introduction}
In non-archimedean functional analysis, there have been various studies on compact-like properties. In general, a non-archimedean valued field is not locally compact, and hence, the classical compactness does not work in non-archimedean functional analysis. Therefore, extracting a part of some properties that are valid for compact sets, we introduce various concepts around the compactness. \par
In section \ref{1a}, using Cauchy filters, we study the c-precompactness. In section \ref{11}, we study the cf-compactness by focusing on the property which every continuous function takes a minimum value on a compact set. As a result, we prove that the c-precompactness is equivalent to the local compactoidity. Furthermore, we prove that the cf-compactness is equivalent to the c-compactness when the whole space is metrizable, and these conditions are not generally equivalent. This result is an answer to the question posed in \cite{varia}. \par
In section \ref{12} and \ref{13}, we study the local compactoidity. The compactoidity is one of the most important concepts in non-archimedean functional analysis, and the local compactoidity is a concept that extends the compacoidity. \par
It is known that a complete local compactoid has good properties (\cite{mod}, \cite{notes}, \cite{local}), but there are few results in the case where the completeness is not imposed. To study a local compactoid in general case, we consider that a local compactoid is roughly the sum of a compactoid and a subspace of finite type. From this perspective, we give some results about non-complete local compactoids in section \ref{12}. Also, we pose some problems which may give light on non-complete local compactoids. In section \ref{13}, we present a result which is already known, but which we are able to prove separately by a non-archimedean functional analysis approach. \par
In sections \ref{gold} and \ref{ebe}, we study the non-archimedean Goldstine theorem and Eberlein-\v{S}mulian theorem, which are related to  compactoidity. The Goldstine theorem for an open unit ball can be proved by the non-archimedean bipolar theorem. On the other hand, the Goldstine theorem for a closed unit ball is unknown. Therefore, we define a Goldstine space, and give a necessary and sufficient condition for a norm-polar space to be a Goldstine space when the coefficient field is spherically complete. As a result, if the valuation group coincides with the set of positive real numbers, then every Goldstine space is finite-dimensional. In section \ref{clo}, the non-archimedean Eberlein-\v{S}mulian theorem is proved, and we apply it to weakly compact operators. \par
If the coefficient field is spherically complete, the closed range theorem was proved in \cite{close}. In section \ref{clo}, we give a new result about the closed range theorem by using epicompactness in the case where the coefficient field is not spherically complete. 

\subsection{Preliminaries}
In this paper, $K$ is a non-archimedean non-trivially valued field which is complete under the metric induced by the valuation $|\cdot| : K \to [0,\infty)$. A unit ball of $K$ is denoted by $B_K := \{ x \in K : |x| \leq 1 \}$.
\par
Throughout,  $E$ is a locally convex space over $K$, that is, a topological vector space over $K$ whose topology is defined by a family of seminorms. A subset $A$ of $E$ is absolutely convex if $A$ is a module over $B_K$. A subset $C$ of $E$ is said to be convex if $C$ is of the form $x + A$ where $A$ is absolutely convex. For a subset $X$ of $E$, using the notation of \cite{lcos}, we denote by $[X]$ the $K$-vector space generated by $X$, by $\mbox{aco} X$ the absolutely convex set generated by $X$, and by $\mbox{co}X$ the smallest convex subset of $E$ containing $X$. \par
Let $A \subseteq E$ be an absolutely convex set. We denote by $A^e$ the edged closure of $A$:
\begin{align*}
  A^e :=  \left\{ 
  \begin{aligned}
    A \quad  &;(\textrm{if $K$ is discretely valued}) \\
    \bigcap_{|\lambda| > 1} \lambda A \quad &; (\textrm{if $K$ is densely valued})
  \end{aligned}
  \right. .
\end{align*}
An absolutely convex set $A$ is called edged if $A^e = A$.
\par
We say that a subset $X$ of $E$ is a compactoid if for every zero neighborhood $U$ in $E$, there exists a finite subset $G \subseteq E$ such that $X \subseteq U + \mbox{aco} G$. An absolutely convex subset $A$ of $E$ is said to be a local compactoid in $E$ if for every zero neighborhood $U$ in $E$, there exists a finite subset $G \subseteq E$ such that $A \subseteq U + [G]$. Note that the local compactoidity depends on an embedding space. On the other hand, if $A$ is a complete local compactoid in $E$, then the local compactoidity of $A$ is independent of an embedding space (see \cite{local}). An absolutely convex set $A$ in a Hausdorff locally convex space $E$ is called c-compact if the following property holds;
\par
Let $\mathcal{C}$ be a family of closed convex subsets of $A$ with the finite intersection property. Then $\bigcap \mathcal{C} \neq \emptyset$.
\par
We write $E'$ for the dual space of $E$, that is, the collection of all continuous linear maps $E \to K$. We denote by $\sigma(E,E')$ the weak topology of $E$ and by $\sigma(E',E)$ the weak star topology of $E'$. We sometimes write $E_{\sigma}'$ instead of $(E',\sigma(E',E))$. Some properties of the weak topology and the weak star topology can be found in \cite{lcos}. Let $F$ be a locally convex space, and $T : E \to F$ be a continuous linear map. The map $f^{*} : F' \to E'$ is defined by $T^{*}(g)(x)= g(T(x))$ for $x \in E$, $g \in F'$. If $(E,\| \cdot \|)$ is a normed space, the norm on $E'$ is defined by the usual operator norm. Moreover, we set $(E'',\|\cdot\|) := (E',\|\cdot\|)'$ and write $j_E$ for the canonical map from $E$ to $E''$.
\par
Let $(E,\| \cdot \|)$ be a normed space. Let $a \in E$, $r > 0$. We write $B_E(a,r)$ for the closed ball with radius $r$ about $a$, that is, $B_E(a,r) := \{ x \in E : \|x-a\| \leq r \}$. A unit ball of $E$ is denoted by $B_E := \{ x \in E : |x| \leq 1 \}$. Let $t \in (0,1]$. A sequence $e_1, e_2, \cdots \in E\backslash \{0\}$ is said to be a $t$-orthogonal sequence (orthogonal sequence for $t = 1$) if for each $n \in \mathbb{N}$, and $\lambda_1, \cdots, \lambda_n \in K$, the inequality
\begin{align*}
t \cdot \max_{1 \leq i \leq n}  \| \lambda_i e_i \| \leq \| \sum_{i=1}^{n} \lambda_i e_i \| \quad 
\end{align*}
holds. \par
A normed space $E$ is of countable type if there exists a subset $X$ that is at most countable such that $\overline{[X]} = E$. In general, a locally convex space $E$ is of countable type if for each continuous seminorm $p$ on $E$, $E_p$ is a normed space of countable type where $E_p$ is the canonical normed space $E/\{p=0\}$ with the norm induced by $p$. We say that $E$ is of finite type if for each continuous seminorm $p$, $E_p$ is a finite-dimensional space.
\par
We say that a seminorm $p$ on $E$ is polar if $p(x)=\sup \{|f(x)| : f \in E' \, \mbox{with} \, |f| \leq p\}$. A locally convex space $E$ is polar if $E$ has a base of polar continuous seminorms, that is, for each continuous seminorm $q$ on $E$, there exists a polar continuous seminorm $p$ which satisfies $q \leq p$. A normed space $(E,\| \cdot \|)$ is norm-polar if the norm $\|\cdot\|$ is a polar seminorm. \par
Let $I$ be an index set. We write $K^I$ for $\prod_I K$, and $K^{(I)}$ for $\{ (x_i)_{i \in I} \in K^I : x_i = 0 \quad \mbox{for all but finitely many} \, i\}$. Let $s : I \to \mathbb{R}_{>0}$ be a map. We define the Banach space $(c_0(I,s),\|\cdot\|)$ by
\begin{align*}
  c_0(I,s) := \{(x_i)_i \in K^I : \lim_i s(i)|x_i| = 0\}, \ \|(x_i)\| := \max_i s(i)|x_i|.
\end{align*}
If $s$ is a constant function whose value is $1$, we write $c_0(I)$ instead of $c_0(I,s)$.

\section{c-precompact} \label{1a}
In this section, we introduce the notion of the c-precompactness which is analogous to the classical precompactness. The main result of this section is Theorem \ref{thmmmm} which shows that the c-precompactness and the local compactoidity coincide.

\begin{dfn} 
 Let $A$ be an absolutely convex set of $E$. A convex filter on $A$ is a filter on $A$ that has a base of convex sets.
\end{dfn}

As in the classical case, we easily obtain the following proposition:

\begin{prop}
 Let $A$ be an absolutely convex set of $E$. Then, $A$ is c-compact if and only if every maximal convex filter on $A$ converges. 
\end{prop} 
As is well known, a metric space $X$ is precompact if and only if every maximal filter on $X$ is a Cauchy filter. Following this observation, we define the c-precompactness. 

\begin{dfn} [Definition of a c-precompact set]
 Let $A$ be an absolutely convex set in $E$. $A$ is said to be a c-precompact set if every maximal convex filter on $A$ is a Cauchy filter.
\end{dfn}

\begin{prop}  [permanence property of c-precompact sets] 
 Let $A$ be a c-precompact set of $E$. Then we have the following:\\
$(1)$ Let $F$ be a locally convex space and $f: E \to F$ be a continuous linear map. Then, $B:=f(A)$ is also c-precompact.\\
$(2)$ The closure $\overline{A}$ of $A$ in $E$ is c-precompact.
\label{propaaa}
\begin{proof}
$(1)$ Let $\mathcal{F}$ be a maximal convex filter on $B$, and $\mathcal{F}'$ a maximal convex filter  on $A$ such that $\mathcal{F}'$ contains $f^{-1}(\mathcal{F}) \cap A := \{ f^{-1}(X) \cap A : X \in \mathcal{F} \}$. Choose an arbitrary absolutely convex zero neighborhood $V$ in $B$, then there exists an absolutely convex zero neighborhood $U$ in $A$  satisfying $f(U)\subset V$. Since $A$ is c-precompact, we can find an element $a$ in $A$  such that $a+U \in \mathcal{F}'$. \par
For any convex set $C\in \mathcal{F}$, we have $(a+U) \cap f^{-1}(C) \neq \emptyset$, so that $(f(a)+f(U)) \cap C \neq \emptyset$. Therefore, we get $f(a)+f(U) \in \mathcal{F}$ by maximality of $\mathcal{F}$. Hence, $f(a) + V \in \mathcal{F}$. Thus, $B$ is c-precompact.
\\
$(2)$ Let $\mathcal{F}$ be a maximal convex filter on $\overline{A}$, and let $U$ be an arbitrary absolutely convex zero neighborhood in $\overline{A}$. Set $\mathcal{F}_0 := \{ (C+U) \cap A : C \in \mathcal{F} \}$. Then there exists a maximal convex filter $\mathcal{F}'$ on $A$ which contains $\mathcal{F}_0$. Since  $A$ is c-precompact,  we can find an element $a$ in $A$  such that $a+(U\cap A)$ is an element of $\mathcal{F}'$. Therefore, for an arbitrary convex set $C \in \mathcal{F}$, we have $(C+U) \cap (a+(U \cap A)) \neq \emptyset$. It follows from the absolute convexity of $U$ that $C \cap (a + U) \neq \emptyset$. Thus, $a+U$ is in $\mathcal{F}$ by maximality of $\mathcal{F}$. Hence, the c-precompactness of $\overline{A}$ follows.  

\end{proof}

\end{prop}

\begin{prop}[{\cite[Lemma 6.3, Theorem 6.7]{notes}}] \label{prop11}
 Let $(E,\| \cdot \|)$ be a normed space, and $A \subset E$ be a complete absolutely convex set. Then the following are equivalent. \\
$(1)$ Each $t$-orthogonal sequence in $A$ tends to $0$. \\
$(2)$ $A$ is a local compactoid in $E$.
\label{propbbb}
\end{prop}

\begin{thm}
 Suppose that $K$ is a spherically complete field. Let $A$ be an absolutely convex set of $E$. Then, the following two conditions are equivalent. \\
$(1)$ $A$ is c-precompact. \\
$(2)$ There exist a locally convex space $G$ and a continuous linear map $f : E \to G$ that is a homeomorphism onto its image such that $f(A)$ is a local compactoid in $G$.
\label{thmmmm}
\begin{proof}
We first prove $(2) \Rightarrow (1)$. Since the c-precompactness does not depend on an embedding space, we may assume that $A$ is a local compactoid in $E$. Let $\mathcal{F}$ be a maximal convex filter on $A$. Choose an arbitrary absolutely convex zero neighborhood $V$ in $E$. Then there exists a finite-dimensional subspace $F \subseteq E$ which satisfies $A \subset V+F$. Since $V \cap F$ is open in $F$ and $F$ is c-compact, $F / V\cap F$ is linearly compact as a $B_K$-module (see \cite[Definition 2.3.10]{mod}). Hence, 
\begin{align*}
  A/V \cap A = (A+V)/V \subset (F+V) / V = F / V\cap F
\end{align*} 
 is also   linearly compact as a $B_K$-module. \par
Consider the maximal convex filter $\mathcal{F} _V := \{ \pi _V (C) : C \in \mathcal{F} \}$ on $A/V \cap A$, where $\pi _ V : A \to A/V \cap A$ is the canonical map. Since $A/V \cap A$ is linearly compact, there exists an element $a \in A$  such that $\{ \pi _V (a) \} \in \mathcal{F} _V$, that is, there exists an element $C \in \mathcal{F}$ such that $\pi _V (C) = \{ \pi _V (a) \}$. Therefore, we have $C \subset a+V \cap A$. Hence $a+V \cap A \in \mathcal{F}$. This means that $\mathcal{F}$ is a Cauchy filter.
\par
Secondly, we prove $(1) \Rightarrow (2)$. Let $\mathcal{P}$ be the set of all continuous seminorms on $E$. For each $p \in \mathcal{P}$, we denote by $\pi _ p : E \to E_{p}^{\land}$ the canonical map where $E_{p}$ is the normed space induced by $p$ and $E_{p}^{\land}$ is the completion of $E_{p}$. From Proposition \ref{propaaa}, we obtain that $\overline{\pi _ p (A)}$ is a c-precompact set.  If $\overline{\pi _ p (A)}$ is a local compactoid in $E_{p}^{\land}$ for each $p \in \mathcal{P}$, it is easily seen that $\prod_{p \in \mathcal{P}} \pi_p(A)$ is a local compactoid in $\prod_{p \in \mathcal{P}} E_{p}^{\land}$. \par
Thus, we may assume that $E$ is a Banach space with a norm $\| \cdot \|$, and $A$ is a closed set in $E$. Hence $A$ is complete. Then by Proposition \ref{propbbb}, it is enough to prove that each $t$-orthogonal sequence in $A$ tends to $0$. \par
Assume that there is a $t$-orthogonal sequence $e_1 , e_2 , \cdots$ in $A$ such that $ \underset{n}{\inf} \| e_n \|  \geq r$ for some $r > 0$ ; we derive a contradiction. Put $C_n := \mbox{co} \{ e_m : m \ge n \}$. Consider a maximal convex filter $\mathcal{F}$ on $A$ which contains all $C_n (n \in \mathbb{N})$. Since $A$ is c-precompact, there exists an element $a \in A$ such that $a+(A \cap B_E(0,tr/2)) \in \mathcal{F}$. In particular, we have $C_1 \cap (a+B_E(0,tr/2)) \neq \emptyset$. By construction of $C_1$ (see \cite[p89]{lcos}), there exist elements $\lambda _1 , \lambda _2 , \cdots , \lambda _ k \in B_K$ such that
\begin{align*}
\lambda _1 + \lambda _2 + \cdots + \lambda _k = 1,  \sum_{1\le i \le k} \lambda_{i}e_i \in a + B_E(0,\dfrac{1}{2}tr).
\end{align*} 
Similarly, it follows from $C_{k+1} \cap (a + B_E(0,tr/2)) \neq \emptyset$ that we can find elements  $\lambda _{k+1} , \lambda _{k+2} , \cdots , \lambda _l  \in B_K (k+1 \leq l)$ satisfying
\begin{align*}
\lambda _{k+1} + \lambda _{k+2} + \cdots + \lambda _l = 1,  \sum_{k+1\le i \le l} \lambda_{i}e_i \in a + B_E(0,\dfrac{1}{2}tr).
\end{align*}
Therefore, we have
\begin{align*}
 \sum_{1\le i \le k} \lambda_{i}e_i - \sum_{k+1\le i \le l} \lambda_{i}e_i \in B_E(0,\dfrac{1}{2}tr).
\end{align*}
On the other hand, since $e_1 , e_2 , \cdots$ is a $t$-orthogonal sequence, we have
\begin{align*}
\|  \sum_{1\le i \le k} \lambda_{i}e_i - \sum_{k+1\le i \le l} \lambda_{i}e_i \|  \geq t\cdot \max_{1\le i \le l}\| \lambda_{i}e_i \| \geq tr\cdot \max_{1\le i \le l} \| \lambda_i \| = tr,
\end{align*}
which is a contradiction.
\end{proof} 

\end{thm}

\begin{cor}[c.f. {\cite[Proposition 2.2]{spc}}]
 Suppose that $K$ is a spherically complete field. Let $A$ be an absolutely convex set of $E$. Then, the following are equivalent. \\
$(1)$ $A$ is c-compact. \\
$(2)$ $A$ is complete and c-precompact. \\
$(3)$ $A$ is a complete local compactoid in $E$.
\label{cornnn}
\end{cor}

\section{cf-compact} \label{11}
In this section, we prove the c-compactness and the cf-compactness (see below) are equivalent in a metrizable space, and are not equivalent in general. This is a problem posed by W. H. Schikhof \cite[p18, Remarks]{varia}. 

\begin{dfn}
 A function $\phi : E \to \mathbb{R}$ is convex if for each $n \in \mathbb{N}$, $x_1, \cdots, x_n \in E$, and $\lambda_1, \cdots, \lambda_n \in B_K$ with $\sum \lambda_i =1$, $\phi$ satisfies the following inequality
\begin{align*}
\phi(\sum_{i=1}^{n}\lambda_ix_i) \leq \max_{1\leq i \leq n}|\lambda_i|\phi(x_i).
\end{align*}
\end{dfn} 

\begin{dfn}
$(1)$ An absolutely convex set $A$ of $E$ is said to be cf-compact if each continuous convex function $\phi : E \to \mathbb{R}$, restricted to $A$, has a minimum.\\
$(2)$ An absolutely convex set $A$ of $E$ is said to be ncf-compact if each continuous non-negative convex function $\phi : E \to [0, \infty)$, restricted to $A$, has a minimum.
\end{dfn}

\begin{prop}
  Let $\phi : E \to \mathbb{R}$ be a convex function. Then, for each $n \in \mathbb{N}$, $x_1, \cdots, x_n \in E$, and $\lambda_1, \cdots, \lambda_n \in B_K$ with $\sum \lambda_i =1$, $\phi$ satisfies the following inequality
\begin{align*}
\phi(\sum_{i=1}^{n}\lambda_ix_i) \leq \max_{1\leq i \leq n} \phi(x_i).
\end{align*}
\end{prop}
\begin{proof}
  Let $n \in \mathbb{N}$, $x_1, \cdots, x_n \in E$, and $\lambda_1, \cdots, \lambda_n \in B_K$ with $\sum \lambda_i =1$. Set
  \begin{align*}
    I := \{ 1 \le n : |\lambda_i| = 1\}, J:= \{1, \cdots, n\} \setminus I.
  \end{align*}
  Then, we have
  \begin{align*}
    \phi(\sum_{i=1}^{n}\lambda_ix_i) &= \phi(\sum_{i \in I} \lambda_i x_i + \sum_{j \in J} (\lambda_j - 1) \cdot x_j + \sum_{j \in J}  x_j) \\
    &\le \max_{i \in I} |\lambda_i|\phi(x_i) \lor \max_{j \in J}|\lambda_j - 1| \phi(x_j) \lor \max_{j \in J} \phi(x_j) \\
    &= \max_{1 \le i \le n} \phi (x_i).
  \end{align*}
\end{proof}

\begin{prop}
  Let $\phi : E \to \mathbb{R}$ be a convex function. Then, for each $t \in \mathbb{R}$, $\phi^{-1}((- \infty, t))$ and $\phi^{-1}((- \infty, t])$ are convex.
\end{prop}
\begin{proof}
  Combine the preceding proposition with \cite[Theorem 3.1.15]{lcos}.
\end{proof}

The following proposition is a simple observation.

\begin{prop}
 A c-compact set is cf-compact. 
\end{prop}

In \cite{varia}, it was proved that in a normed space, the c-compactness is equivalent to the cf-compactness. The proof of \cite{varia} shows that a stronger claim holds:

\begin{prop}[{\cite[Theorem 4.2]{varia}}]
 Let $A$ be an absolutely convex set of a normed space $E$. Then, the following are equivalent.\\
$(1)$ $A$ is c-compact. \\
$(2)$ $A$ is cf-compact. \\
$(3)$ $A$ is ncf-compact.
\label{propeee}
\end{prop}

To examine the cf-compactness, we study some properties of convex functions.

\begin{prop}
 Let $\phi : E \to \mathbb{R}$ be a convex function. Then, for each $x_0 \in E$, $\psi(x):=\phi(x+x_0)$ is also a convex function.
\end{prop}

\begin{prop}
 Let $\psi : [0,\infty) \to (-\infty,0]$ be a monotone increasing function. Then, for each seminorm $p$ on $E$,
\begin{align*}
\phi : E \ni x \mapsto \psi(p(x)) \in \mathbb{R}
\end{align*}{}
is a convex function.
\label{propzzz}
\begin{proof}
For $n \in \mathbb{N}$, $x_1, \cdots, x_n \in E$, $\lambda_1, \cdots, \lambda_n \in B_K$ with $\sum \lambda_i =1$, we have
\begin{align*}
\begin{split}
\phi(\sum_i \lambda_i x_i) = \psi(p(\sum_i \lambda_i x_i)) &\leq \psi(\max_i |\lambda_i|\cdot p(x_i))
\\
&= \max_i \psi(|\lambda_i|\cdot p(x_i))
\\
&\leq \max_i \psi(p(x_i))
\\
&\leq \max_i |\lambda_i|\cdot \psi(p(x_i)) = \max_i |\lambda_i|\cdot \phi(x_i).
\end{split}
\end{align*}{}
Notice that the last inequality holds since $\psi$ is a non-positive function.
\end{proof}
\end{prop}

Now, we are ready to prove the main theorem of this section which expands Proposition \ref{propeee}.

\begin{thm}
 Let $A$ be an absolutely convex set of a metrizable locally convex space $E$. Then, the following are equivalent. \\
$(1)$ $A$ is cf-compact.
\\
$(2)$ $A$ is c-compact.
\label{thmfff}
\begin{proof}
It suffices to prove $(1) \Rightarrow (2)$. For each $n \in \mathbb{N}$, we define a monotone increasing function $\psi_n$ on $[0,\infty)$ as follows:
\begin{align*}
\psi_n(x) := -1 + \frac{1}{n} (1-e^{-x}).
\end{align*}{}
Notice that $\psi_n$ is a non-positive function. Secondly, let $p_1\leq p_2 \leq \cdots$ be continuous seminorms defining the topology of $E$ and set $\phi(x) := \max_n \psi_n(p_n(x))$ $(x\in E)$. From Proposition \ref{propzzz}, $\phi$ is a convex function and takes the minimum value $-1$ only at $0$. It is easily seen that $\phi$ is continuous . 
Using this continuous convex function $\phi$, we obtain the desired result as follows.
\par
Let $A$ be a cf-compact set. Clearly we may assume $A$ is not a singleton. From Proposition \ref{propeee}, for each continuous seminorm $p$, $\pi_p(A)$ is c-compact where $\pi_p : E \to E_p$ is the canonical map. Hence, $K$ is spherically complete. From this observation, $\pi_p(A)$ is a local compactoid in $E_p$ for each continuous seminorm $p$, and therefore, $A$ is a local compactoid in $E$. \par 
Thus, from Corollary \ref{cornnn}, it is enough to prove that $A$ is complete. Moreover, by taking the completion of $E$, it suffices to prove that $A$ is a closed set. Let $x_0$ be an element of $\overline{A}$. Then, $\psi(x) := \phi(x-x_0)$ is a continuous convex function on $E$. Since $A$ is cf-compact, we have the equality 
\begin{align*}
\min \psi(A) = \inf \psi(A) = \inf \psi(\overline{A})= \psi(x_0)
\end{align*}
This equality implies $x_0 \in A$ because $\psi$ takes its minimum value only at $x_0$ by a choice of $\phi$.
\end{proof}
\end{thm}

\begin{prop}
 Let $\phi : K \to [0,\infty)$ be a convex function taking the minimum value at $0$ and assume that $\phi$ is not constant. Then, $\lim_{|x|\to \infty}\phi(x)= \infty$.
\label{propggg}
\begin{proof}
Let $x,y \in K$ with $|x| \le |y|$. Since $\phi$ is a convex function, $\phi^{-1}((-\infty,\phi(y)])$ is an absolutely convex set containing $y$. Thus, we have $x \in \phi^{-1}((-\infty,\phi(y)])$, that is, $\phi(x) \le \phi(y)$. Therefore, $M:=\lim_{|x|\to \infty}\phi(x)$ exists. \par
Suppose $M<\infty$, we derive a contradiction. By assumption, there exists an element $x \in K$ such that $\phi(0) < \phi(x)$. Let $\lambda \in B_K$ with $\lambda \neq 0$. Then it follows that
\begin{align*}
\phi(x) = \phi(\lambda \cdot \frac{x}{\lambda} + (1-\lambda)\cdot 0) \leq |\lambda|\phi(\frac{x}{\lambda}) \leq |\lambda|M \to 0 \quad \mbox{as} \quad \lambda \to 0,
\end{align*}
which is a contradiction.
\end{proof}
\end{prop}

The following example claims that $E$ must be a normed space in Proposition \ref{propeee}.

\begin{exam} \label{exhhh}
Let $K$ be spherically complete, and let $E := K^{\mathbb{N}}$ equipped with the product topology. Then, there exists an ncf-compact set $A$ of $E$ which is not cf-compact.
\begin{proof} 
Let $A := K^{(\mathbb{N})} \subseteq E$. This absolutely convex set $A$ satisfies the desired property. Indeed, since $A$ is not closed, $A$ is not c-compact. Hence, it follows from Theorem \ref{thmfff} that $A$ is not cf-compact. \par 
It remains to prove that $A$ is ncf-compact. Let $\phi$ be a non-negative continuous convex function on $E$. We have to prove that $\min \phi(A)$ exists. By the c-compactness of $E$, there exists an element $x_0 \in E$ such that $\phi$ attains its minimum at $x_0$. Let $\psi : E \to [0,\infty)$ be a continuous convex function defined by the formula $\psi(x) := \phi(x+x_0)$. Since $\psi$ is continuous, there exists a sufficiently large number $N \in \mathbb{N}$ such that $\psi(\{ (x_n)_{n \in \mathbb{N}} \in E : x_n = 0 \; \mbox{for all} \;  1 \leq n \leq N\})$ is bounded. From Proposition \ref{propggg}, we obtain $\psi(\{ (x_n)_{n \in \mathbb{N}} \in E : x_n = 0 \; \mbox{for all} \;  1 \leq n \leq N\}) = \{\psi(0)\}$. Now, it is clear that $\min \phi(A)$ exists.
\end{proof}
\end{exam}

The following example shows that the metrizablity condition of $E$ can not be dropped in Theorem \ref{thmfff}.

\begin{exam}
Let $K$ be spherically complete, $I$ be an uncountable set and,  let $E := K^{I}$ equipped with the product topology. Then, there exists a cf-compact set $A$ of $E$ which is not c-compact.
\begin{proof}
We will see that 
\begin{align*}
A := \{ (x_i)_{i \in I} \in E : \{i \in I : x_i \neq 0\} \, \mbox{is at most countable}\}
\end{align*}
is the desired example. It is clear that $A$ is not c-compact. We prove that $A$ is cf-compact. Let $\phi$ be a continuous convex function on $E$. By an argument similar to that as in Example \ref{exhhh}, we define a convex function $\psi(x):=\phi(x+x_0)$ where $x_0 \in E$ is the element which minimizes the function $\phi$. For each $n \in \mathbb{N}$, there exists a finite subset $J_n$ of $I$ satisfying 
\begin{align*}
\psi(\{(x_i)_{i \in I} \in E : x_i = 0 \quad \mbox{for all}  \quad i \in J_n\}) \subseteq (-\infty , \psi(0) + \frac{1}{n}).
\end{align*}{}
Set $J := \bigcup J_n$. Then $J$ is countable and satisfies the following equality
\begin{align*}
\psi(\{(x_i)_{i \in I} \in E : x_i = 0 \quad \mbox{for all}  \quad i \in J\}) = \{\psi(0)\}. 
\end{align*}{}
Thus, it is easy to see that $\min \phi(A)$ exists.
\end{proof}
\end{exam}

\section{some problems of local compactoids} \label{12}

Due to \cite{local}, a complete local compactoid is the sum of a compactoid and  a subspace of finite type. In addition, completeness condition can not be dropped. On the other hand, we have the following theorem. 

\begin{thm}
 Let $(E , \| \cdot \|)$ be a normed space and $A$ be a local compactoid in $E$. Then, there exist a finite dimensional subspace $F \subseteq E$ and a compactoid $B \subseteq E$ which satisfy $A \subseteq F+B$. 
\begin{proof}
By the local compactoidity of $A$, for each $n \in \mathbb{N}$, there exists a finite dimensional subspace $F_n \subseteq E$ with $A \subseteq F_n +B_{E}(0,1/n)$. Put $E_0 := [A]+[\bigcup_n F_n]$. Since $A$ is a local compactoid, $E_0$ is a normed space of countable type. Therefore, $F_1$ is complemented in $E_0$, so that we can choose a subspace $E_1 \subseteq E_0$ which is topologically complemented to $F_1$. \par
  Let $P_1 : E_0 \rightarrow F_1$ and $P_2 : E_0 \rightarrow E_1$ be the canonical continuous projections. Take an element $a \in A$. Since $A$ is contained in $F_1 +B_{E_0}(0,1)$, there exist elements $f \in F_1$ and $v \in B_{E_0}(0,1)$ satisfying $a = f+v$. Then, we have $\| P_2 (a) \| = \| P_2 (v) \| \leq \| P_2 \|$, where $\| P_2 \|$ is the operator norm of $P_2$. Therefore, $P_2 (A)$ is a local compactoid in $E_1$ and bounded. Hence, $P_2 (A)$ is a compactoid (c.f. \cite[Theorem 3.8.4]{lcos}). Finally, we have $A = P_1 (A) + P_2 (A) \subseteq F_1 + P_2 (A)$.    
\end{proof}
\end{thm}

\begin{prob*}
Let $A$ be a local compactoid in $E$ (not necessary a normed space). Then, are there a subspace $F \subseteq E$ of finite type and a compactoid $B \subseteq E$ such that $A \subseteq F + B$ ? 
\end{prob*}

Note that the converse is clearly true. Therefore, the above condition will be a new characterization of the local compactoidity if the problem is solved. \par 
\vspace{5pt}
Let $A$ be an absolutely convex set of $E$. As we have already noted, the local compactoidity depends on an embedding space. On the other hand, if $A$ is the sum of a compactoid and a space of finite type, the local compactoidity of $A$ is independent of an embedding space. We will study the sufficient condition for ``$A$ is the sum of a compactoid and a space of finite type''.  


\begin{exam}[compare {\cite[Example 3.6]{local}}]
 There exist a normed space $(E , \| \cdot \|)$ and a local compactoid $A$ in $E$ which satisfies the following two conditions: \\
$(1)$ $A$ is not the sum of a compactoid and a space of finite type. \\
$(2)$ Let $F$ be a locally convex space, and let $B$ be an absolutely convex set contained in $F$. Suppose that there exists an isomorphism $\phi : A \to B$ as topological $B_K$-modules. Then, $B$ is also a local compactoid in $F$. 
\begin{proof}
Let $E := (K^{(\mathbb{N})},\| \cdot \|)$ with the norm defined by the formula $\|(x_{n})_{n \in \mathbb{N}}\| := \max |x_n|$. We write $e_1, e_2, \cdots$ for the canonical unit vectors. Put $f_n := e_{2^n}, g_{(i,j)} := e_{3^i5^j}$ and
\begin{align*}
A := \text{aco} \{e_1\} + \text{aco} \{ \lambda^n e_1 + \lambda^{-n} f_n : n \in \mathbb{N} \} + \text{aco} \{\lambda^{-n}f_{n} + \lambda^{-j} g_{(n,j)} : n, j \in \mathbb{N}\}
\end{align*}
where $\lambda \in K$ with $|\lambda| > 1$. We prove that $A$ satisfies the desired conditions. \par
  Suppose that $A$ is the sum of a compactoid and a space of finite type. Let $P_n : E \to K$ be the canonical $n$-th coordinate projection. Then, $P_n(A)$ is bounded for each $n \neq 1$. Therefore, since $A$ is unbounded, $A$ must contain $[e_1]$, which is a contradiction and we obtain $(1)$. \par
  Secondly, we prove $(2)$. Let $F, B$, and $\phi$ be as in $(2)$. For each $n,j \in \mathbb{N}$, we have  $\lambda^{-2n}f_n, \lambda^{-n-j} g_{(n,j)} \in A$. Since $F$ is a locally convex space, we obtain
\begin{align*}
\phi(\lambda^{-n}f_{n} + \lambda^{-j}g_{(n,j)}) = \lambda^{n}\{\phi(\lambda^{-2n}f_{n}) + \phi(\lambda^{-n-j}g_{(n,j)})\} \to \lambda^{n}\phi(\lambda^{-2n}f_n) \ \text{as} \ j \to \infty.
\end{align*}
Let $C := \text{aco} \{\phi(\lambda^{-n}f_{n} + \lambda^{-n-j}g_{(n,j)}) : n, j \in \mathbb{N}\}$, then $C$ is a compactoid and the closure $\overline{C}$ of $C$ in $F$ contains $\lambda^{n}\phi(\lambda^{-2n}f_n)$ for each $n \in \mathbb{N}$. Finally, since $\phi(\lambda^n e_1 + \lambda^{-n} f_n) = \lambda^n\phi(e_1) + \lambda^{n}\phi(\lambda^{-2n}f_n)$ for each $n \in \mathbb{N}$, we have $B \subseteq [\phi(e_1)] + \overline{C}$, which implies $(2)$. 
\end{proof}
\end{exam}

\begin{thm}
 Let $(E , \| \cdot \|)$ be a normed space and $A$ be a local compactoid in $E$. Then $[A]$ is not a Banach space whenever $[A]$ is an infinite dimensional space.
\begin{proof}
Suppose that $[A]$ is a Banach space. Then, $A$ is a local compactoid in $[A]$. Therefore, there exist elements $a_1, \cdots , a_n \in A$, and a compactoid $B \subseteq [A]$ satisfying $A \subseteq [a_1, \cdots , a_n] + B$. Since $[A]$ is a normed space of countable type, there exists a subspace $F \subseteq [A]$ which is topologically complemented to $[a_1, \cdots , a_n]$ in $[A]$. Let $P : [A] \rightarrow F$ be the continuous projection. Then, we have $F = [P (B)]$. As $F$ is a Banach space, the closure of $P (B)$ is a compactoid zero neighborhood in $F$. Hence, $F$ is a finite-dimensional space by the Riesz lemma (\cite[Theorem 3.8.6]{lcos}).   
\end{proof}
\end{thm}

\begin{cor}
 Let $(E , \| \cdot \|)$ be a normed space and $A$ be a local compactoid in $E$. If $[A]$ is a Banach space, then $A$ is the sum of a compactoid and a space of finite type.
\end{cor}

\section{topological stability of complete metrizable local compactoids} \label{13}

In this section, we will prove Theorem \ref{thmooo}. This result has been already proved in \cite[Theorem 5.4.8]{mod}. We will give another proof which is based on the weak topology. It has been known that complete metrizable compactoids and c-comact sets are topologically stable (Propositions \ref{propjjj} and \ref{propppp}).

\begin{prop}[{\cite[Theorem 3.2]{top}}]
 Let $A$ be a complete absolutely convex metrizable compactoid in a Hausdorff locally convex space $(E,\tau)$. Let $\tau_1$ be a Hausdorff locally convex topology on $E$, weaker than or equal to $\tau$. Then $\tau = \tau_1$ on $A$. 
\label{propjjj}
\end{prop}

\begin{prop}[{\cite[Proposition 3.2]{spc}}]
 Let $A$ be a c-compact set in a Hausdorff locally convex space $(E,\tau)$. Let $\tau_1$ be a Hausdorff locally convex topology on $E$, weaker than or equal to $\tau$. Then $\tau = \tau_1$ on $A$.
\label{propppp}
\end{prop}

\begin{cor}
 Let $I$ be an index set and $(K^I, \tau )$ be a locally convex space equipped with the product topology. Let $\tau_1$ be a Hausdorff locally convex topology on $K^I$, weaker than or equal to $\tau$. Assume that $K$ is spherically complete. Then, we have $\tau = \tau_1$. 
\end{cor}

As the above proposition, we have the following theorem.

\begin{thm}
 Let $(K^{\mathbb{N}}, \tau )$ be a locally convex space equipped with the product topology. Let $\tau_1$ be a Hausdorff locally convex topology on $K^{\mathbb{N}}$, weaker than or equal to $\tau$. Then, we have $\tau = \tau_1$. 
\label{thmiii}
\begin{proof}
Let $K^{(\mathbb{N})}$ be a locally convex space equipped with the strongest locally convex topology. Let $F := (K^{\mathbb{N}}, \tau_1 )' \subseteq (K^{\mathbb{N}}, \tau )'=K^{(\mathbb{N})}$. Since $(K^{\mathbb{N}}, \tau_1 )$ is a Hausdorff locally convex space of finite type, $F$ is weakly dense in $K^{(\mathbb{N})}$. As $K^{(\mathbb{N})}$ is of countable type and $F$ is a subspace of $K^{(\mathbb{N})}$, the weak closure of $F$ and the closure of $F$ coincide. In the strongest locally convex topology, every convex set is closed. Hence, $F$ is equal to $K^{(\mathbb{N})}$, and we obtain that $\tau = \tau_1$. 
\end{proof}
\end{thm}

\begin{prop}[{\cite[Proposition 1.6]{local}}]
 Let $A$ be a local compactoid in a polar locally convex space $E$. Then, on $A$, the weak topology $\sigma (E, E')$ and the initial topology coincide.
\label{propkkk}
\end{prop}

\noindent \textbf{Construction} \par
Let $(E,\tau)$ be a polar Hausdorff locally convex space. Set 
\begin{align*}
\gamma_{\tau} : E \longrightarrow (K^{((E,\tau)')})_{\sigma}' \, ; \, \gamma_{\tau}(x)((a_f)_{f \in (E,\tau)'}):= \sum a_f f(x) \quad (x \in E).
\end{align*}
where $K^{((E,\tau)')}$ is equipped with the strongest locally convex topology. Obviously, the right-hand side is a finite sum. From Proposition \ref{propkkk}, if $A$ is a local compactoid in $E$, the restriction $\gamma_{\tau}|_A$ of $\gamma_{\tau}$ to $A$ is a homeomorphism onto its image. \\
Let $I$ be an index set and $K^{(I)}$ be a locally convex space equipped with the strongest locally convex topology. For $X \subseteq K^{(I)}$, $Y \subseteq (K^{(I)})'$ we set
\begin{align*}
\begin{split}
X^{^{\circ}} &:= \{f \in (K^{(I)})' : |f(x)| \leq 1 \quad \mbox{for all} \quad x \in X\} \\
Y^{^{\circ}} &:= \{x \in K^{(I)} : |f(x)| \leq 1 \quad \mbox{for all} \quad f \in Y\}.
\end{split}
\end{align*}{}
As is well known, each closed subspace $B \subseteq (K^{(I)})_{\sigma}'$ satisfies the equality $B^{^{\circ \circ}} = B$. \par
\vspace{4pt}
We consider the following situation:\\
Let $A$ be a local compactoid in a polar Hausdorff locally convex space $(E,\tau)$, and $D$ be a complete subspace of $A$. We fix an algebraic complement $F \subseteq K^{((E,\tau)')}=:H$ of $\gamma_{\tau}(D)^{^{\circ}}$, that is, $H=F \oplus \gamma_{\tau}(D)^{^{\circ}}$ algebraically. Let
\begin{align*}
\pi_1 : H \longrightarrow F \subseteq H, \quad \pi_2 : H \longrightarrow \gamma_{\tau}(D)^{^{\circ}} \subseteq H
\end{align*}{}
be the canonical continuous projections.

\begin{prop}[c.f. {\cite[Proposition 2.3]{local}}]
 Using the above notations, we have the following: \\
$(1)$ $\pi_1^{*}(\gamma_{\tau}(A))$, $\pi_2^{*}(\gamma_{\tau}(A)) \subseteq \gamma_{\tau}(A)$. \\
$(2)$ $\pi_1^{*} \times \pi_2^{*} : \gamma_{\tau}(A) \to \pi_1^{*}(\gamma_{\tau}(A)) \times \pi_2^{*}(\gamma_{\tau}(A))$ is a homeomorphism, where the right-hand side is equipped with the product topology.\\
$(3)$ $\textup{Ker}\,\pi_2^{*}|_{\gamma_{\tau}(A)}= \pi_1^{*}(\gamma_{\tau}(A)) = \gamma_{\tau}(D)$. \\
Moreover, suppose that $A$ is complete and $D$ is the largest subspace of $A$. Then, $D$ is closed and we have \\
$(4)$ $\textup{Ker}\,\pi_1^{*}|_{\gamma_{\tau}(A)}$ is a complete compactoid.
\label{proplll}
\begin{proof}
$(1)$ Since $\pi_1^{*} + \pi_2^{*} = \mbox{id}_{H'}$, it is enough to prove that $\pi_1^{*}(\gamma_{\tau}(A))$ is contained in $\gamma_{\tau}(A)$. Let $f \in \gamma_{\tau}(A)$. Then we have $\pi_1^{*}(f)(\gamma_{\tau}(D)^{^{\circ}})=0$. Therefore it follows that $\pi_1^{*}(f) \in \gamma_{\tau}(D)^{^{\circ \circ}} =\gamma_{\tau}(D)$. \\
$(2),(3)$ It is easy to prove. \\
$(4)$ See \cite[Proposition 2.3]{local}.
\end{proof}
\end{prop}

\begin{thm}
 Let $A$ be a complete metrizable local compactoid in a Hausdorff locally convex space $(E,\tau)$. Let $\tau_1$ be a Hausdorff locally convex topology on $E$, weaker than or equal to $\tau$. Then $\tau = \tau_1$ on $A$.
\label{thmooo}
\begin{proof}
Since $A$ is complete, we may assume $E = [A]$. In particular, $E$ is a space of countable type. Hence, $E$ is a polar space. We can consider $(E,\tau_1)' \subseteq (E,\tau)'$ and let $\iota : H_1:=K^{((E,\tau_1)')} \to H:=K^{((E,\tau)')}$ be the canonical injection. Then, we have the following commutative diagram
\[
\xymatrix@M=8pt{
(E,\tau) \ar[r]^-{\gamma_{\tau}} \ar[d]^-{\mbox{id}} & H' \ar[d]^-{\iota^{*}} &\\
(E,\tau_1)\ar[r]^-{\gamma_{\tau_{1}}} & H_1' & .
}
\]
Let $D$ be the largest subspace of $A$. Then, $D$ is a complete metrizable local compactoid space with respect to $\tau$. Therefore, $D$ is a complete space of finite type. Hence, $D$ is of the form $K^{\mathbb{N}}$ with the product topology (\cite[Remark 5.3.12]{lcos}). From Theorem \ref{thmiii}, $\tau$ and $\tau_1$ coincide on $D$. In particular $D$ is also complete with respect to $\tau_1$. \par
By construction, it is easy to see that $\iota(\gamma_{\tau_{1}}(D)^{^{\circ}}) = \iota(H_1) \cap \gamma_{\tau}(D)^{^{\circ}}$.  Let $F_1 \subseteq H_1$ be an algebraic complement of $\gamma_{\tau_{1}}(D)^{^{\circ}}$. It follows from $\iota(\gamma_{\tau_{1}}(D)^{^{\circ}}) = \iota(H_1) \cap \gamma_{\tau}(D)^{^{\circ}}$ that we can choose $F \subseteq H$ an algebraic complement of $\gamma_{\tau}(D)^{^{\circ}}$ containing $\iota(F_1)$. Let
\begin{align*}
\begin{split}
&\pi_1 : H \longrightarrow F \subseteq H, \quad \pi_2 : H \longrightarrow \gamma_{\tau}(D)^{^{\circ}} \subseteq H, \\
&p_1 : H_1 \longrightarrow F_1 \subseteq H_1, \quad p_2 : H_1 \longrightarrow \gamma_{\tau_1}(D)^{^{\circ}} \subseteq H_1
\end{split}
\end{align*}{}
be the canonical projections. Then the following diagram is commutative
\[
\xymatrix@M=15pt{
H_1 \ar[r]^-{p_1 \times p_2} \ar[d]^-{\iota} & F_1 \times \gamma_{\tau_1}(D)^{^{\circ}} \ar@<-0.5ex>[d]^-{\iota \times \iota} & \\
 H \ar[r]^-{\pi_1 \times \pi_2} & F \times \gamma_{\tau}(D)^{^{\circ}} & .
}
\]
Therefore, we have the following commutative diagram
\[
\xymatrix@M=15pt{
\gamma_{\tau}(A) \ar[r]^-{\pi_1^{*} \times \pi_2^{*}} \ar[d]^-{\iota^{*}} & \pi_1^{*}(\gamma_{\tau}(A)) \times \pi_2^{*}(\gamma_{\tau}(A)) \ar[d]^-{\iota^{*} \times \iota^{*}} & \\
\gamma_{\tau_1}(A) \ar[r]^-{p_1^{*} \times p_2^{*}} & p_1^{*}(\gamma_{\tau_1}(A)) \times p_2^{*}(\gamma_{\tau_1}(A)) & .
}
\]
From Proposition \ref{proplll}, both horizontal arrows are homeomorphisms. Thus, to prove the theorem, it is enough to prove that the right vertical arrow is a homeomorphism. By Proposition \ref{proplll}, we have $\pi_1^{*}(\gamma_{\tau}(A))=\gamma_{\tau}(D)$, and $p_1^{*}(\gamma_{\tau_1}(A))=\gamma_{\tau_1}(D)$. Hence, by Theorem \ref{thmiii}, $\iota^{*} : \pi_1^{*}(\gamma_{\tau}(A)) \to p_1^{*}(\gamma_{\tau_1}(A))$ is a homeomorphism. Finally, let $B:= \mbox{Ker}\,(\pi_1^{*}|_{\gamma_{\tau}(A)} \circ \gamma_{\tau}|_A)$. Then we get $\pi_2^{*}(\gamma_{\tau}(A))=\gamma_{\tau}(B)$, and $p_2^{*}(\gamma_{\tau_1}(A))=\gamma_{\tau_1}(B)$. From Propositions \ref{propjjj} and \ref{proplll}, $B$ is a complete metrizable compactoid, and $\iota^{*} : \pi_2^{*}(\gamma_{\tau}(A)) \to p_2^{*}(\gamma_{\tau_1}(A))$ is a homeomorphism. Therefore, the proof is complete.
\end{proof}
\end{thm}

\section{The Goldstine theorem} \label{gold}

In this section, we study the Goldstine theorem. The non-archimedean Goldstine theorem for an open unit ball is known.

\begin{thm}[c.f. {\cite[Proposition 3.4]{equ}}] \label{ThmG1}
 Let $(E,\|\cdot\|)$ be a norm-polar Banach space. Then, we have
 \begin{align*}
   \left(\overline{j_E(B_E)}^{\sigma(E'',E')}\right)^{e} = B_{E''}.
 \end{align*}
\end{thm}

\begin{rem}
  Even if $K$ is discretely valued, the norm-polarity of $E$ is necessary. Indeed, let $E = c_0 (\mathbb{N},s)$ for which $s(1)>s(2)> \cdots \to 1$ as $n \to \infty$. Let $f \in (c_0 (\mathbb{N},s))'$ with $f(e_n) = 1$ for each $n \in \mathbb{N}$ where $e_1,e_2, \cdots$ are the canonical unit vectors. Then, it is easy to see $\|f\| = 1$. Thus, there exists $\theta \in B_{E''}$ such that $\theta(f) = 1$. \par
On the other hand, for each $x = (x_1,x_2,\cdots) \in B_E$, we have $|x_n| \le |\pi|$ for each $n \in \mathbb{N}$ where $\pi \in B_K$ is a generating element of a maximal ideal of $B_K$. Therefore, we get $j_E(x)(f) \le |\pi|$ for each $x \in B_E$. This implies $\theta \notin \left(\overline{j_E(B_E)}^{\sigma(E'',E')}\right)^{e}$.
\end{rem}

We shall study the Goldstine theorem for a closed unit ball. 

\begin{dfn}
  Let $(E,\|\cdot\|)$ be a norm-polar Banach space. Then, $E$ is said to be a Goldstine space if 
  \begin{align*}
    \overline{j_E(B_E)}^{\sigma(E'',E')} = B_{E''}.
  \end{align*}
\end{dfn}

\begin{rem}
  Reflexive Banach spaces are Goldstine spaces. In particular, if $K$ is not spherically complete, then $(c_0,\|\cdot\|)$ is a Goldstine space.
\end{rem}

Since $j_E(B_E)$ is an edged set, by Theorem \ref{ThmG1}, it is natural to ask the following question (\cite[p95, Problem]{lcos}):
\begin{center}
  Does it follow that the closure of an edged set is also edged?
\end{center}
We have a negative answer to it.

\begin{exam}
  Let $\omega$ be a symbol, and let $(E,\|\cdot\|) = c_0((\mathbb{N}\times \mathbb{N}) \amalg \{\omega\})$. Take a sequence $a_1, a_2, \cdots \in K$ such that $|a_1| < |a_2| < \cdots \to 1$ as $n \to \infty$ and $\rho \in K$ with $0 < |\rho| < 1$. Let $e_{\omega},e_{(i,j)}$ $i,j \in \mathbb{N}$, be canonical unit vectors. Now, let us put
  \begin{align*}
    f_{(i,j)} := a_i e_{\omega} + \rho^j e_{(i,j)} \ (i,j \in \mathbb{N})
  \end{align*}
  and
  \begin{align*}
    A := \mathrm{aco}\{f_{(i,j)} : i,j \in \mathbb{N}\}.
  \end{align*}
  Then, $A$ is an edged set, but $\overline{A}$ is not edged.
\end{exam}
\begin{proof}
  Since the set $\{f_{(i,j)} : i,j \in \mathbb{N}\}$ is linearly independent over $K$, it is clear that $A$ is an edged set. Let us prove that $\overline{A}$ is not edged. For each $i \in \mathbb{N}$, $\overline{A}$ contains $a_i e_{\omega}$. Thus, $(\overline{A})^e$ contains $e_{\omega}$. On the other hand, we see $\|e_{\omega} - a\| \ge 1$ for each $a \in A$. Hence, we have $e_{\omega} \notin \overline{A}$, which completes the proof.
\end{proof}

By the above example, it is not trivial whether a norm-polar Banach space is a Goldstine space or not.

\begin{thm} \label{ThmG2}
  Let $K$ be spherically complete and densely valued. Then, for a Banach space $(E,\|\cdot\|)$, the following are equivalent: \\
  $(1)$ $E$ is a Goldstine space \\
  $(2)$ For each $f \in E'$, $f(B_E)$ is an edged set.
\end{thm}
\begin{proof}
  First, let us suppose $(1)$. Then, $\overline{j_E(B_E)}^{\sigma(E'',E')} = B_{E''}$ is an edged set. Moreover, by $p$-adic Alaoglu theorem (c.f. \cite[Proposition 3.1]{equ}), we have that $\overline{j_E(B_E)}^{\sigma(E'',E')}$ is c-compact. Thus, by \cite[Theorem 6.2.3]{lcos}, 
  \begin{align*}
    \{\theta(f) : \theta \in \overline{j_E(B_E)}^{\sigma(E'',E')}\}
  \end{align*}
is edged for each $f \in B_E$. Since each convex set in $K$ is closed, we get 
\begin{align*}
  \{\theta(f) : \theta \in \overline{j_E(B_E)}^{\sigma(E'',E')}\} \subseteq \overline{\{\theta(f) : \theta \in j_E(B_E)\}} = f(B_E).
\end{align*}
Therefore, we have that $f(B_E)$ is an edged set. \par
Conversely, suppose that there exists $\psi \in B_{E''}$ for which 
\begin{align*}
  \psi \notin \overline{j_E(B_E)}^{\sigma(E'',E')}.
\end{align*}
Then, there exists $f \in E'$ such that $\psi(f) = 1$ and $|\theta (f)| < 1$ for each $\theta \in \overline{j_E(B_E)}^{\sigma(E'',E')}$. Therefore, it follows from Theorem \ref{ThmG1} that we have
\begin{align*}
  f(B_E) = \{\theta(f) : \theta \in \overline{j_E(B_E)}^{\sigma(E'',E')}\} = \{x \in K : |x| < 1\},
\end{align*}
which completes the proof.
\end{proof}

As an application of the above theorem, we have the following.

\begin{thm}
  Let $K$ be spherically complete, and suppose $|K| = \mathbb{R}_{\ge 0}$. Then, every Goldstine space is finite-dimensional.
\end{thm}
\begin{proof}
  Let $(E,\|\cdot\|)$ be a Goldstine space. Then by Theorem \ref{ThmG2}, $f(B_E)$ is an edged set for each $f \in E'$. From the assumption $|K| = \mathbb{R}_{\ge 0}$, we have that $B_E$ is a weakly c'-compact set (see \cite{wc}). Then by \cite[Theorem 5.2.13]{wc}, $B_E$ is a compactoid. Thus, by the Riesz lemma, $E$ is finite-dimensional.
\end{proof}

\begin{prob*}
  Let $K$ be spherically complete and densely valued. Suppose $|K| \neq \mathbb{R}_{\ge 0}$. Then, does it follow that there exists an infinite-dimensional Goldstine space?
\end{prob*}

The above problem is open, while we have the following.

\begin{prop}
  Let $K$ be spherically complete and densely valued. Let $I$ be an index set with $\# I = \infty$ and $s : I \to \mathbb{R}_{>0}$. Then, $c_0(I,s)$ is not a Goldstine space.
\end{prop}
\begin{proof}
  Let $\lambda \in K$ with $|\lambda| > 1$. To prove the theorem, we may assume $s(I) \subseteq [1,|\lambda|]$. Then, there exist distinct elements $i_1, i_2, \cdots \in I$ and $r \in [1,|\lambda|]$ such that $s(i_n) \to r$ as $n \to \infty$. For each $n \in \mathbb{N}$, choose $\mu_n \in K$ such that
  \begin{align*}
    1 - \frac{1}{n} < |\mu_n| \cdot \frac{1}{s(i_n)} < 1.
  \end{align*}
  Let us define $f \in E'$ by
  \begin{align*}
    f(e_i) := \left\{
  \begin{aligned}
    \mu_n \quad &;(i = i_n \ \mathrm{for} \ \mathrm{some} \ n \in \mathbb{N}) \\
    0 \quad &; \mathrm{otherwise}
  \end{aligned}
  \right. ,
  \end{align*}
  where $e_i$'s, $i \in I$, are the canonical unit vectors. Then, we have
  \begin{align*}
    f(B_{c_0(I,s)}) = \{x \in K : |x| < 1\}.
  \end{align*}
  Thus, by Theorem \ref{ThmG2}, $c_0(I,s)$ is not a Goldstine space.
\end{proof}

\section{The Eberlein-\v{S}mulian theorem} \label{ebe}

The non-archimedean Eberlein-\v{S}mulian theorem for the weak compactness was studied in \cite{eber}. On the other hand, we shall study the weakly complete compactoids.

\begin{thm} \label{thme1}
  Let $K$ be spherically complete, and $A$ be an absolutely convex subset of a Banach space $(E,\|\cdot\|)$. Suppose that $A$ is weakly complete. Then, $A$ is a norm complete local compactoid.
\end{thm}
\begin{proof}
  By \cite[Theorem 5.5.2]{lcos}, $A$ is norm complete. Suppose that $A$ is not a local compactoid. Then by Proposition \ref{prop11}, there exist $t \in (0,1)$ and a bounded $t$-orthogonal sequence $f_1,f_2, \cdots \in A$ such that $\inf_n \|f_n\| > 0$. Put 
  \begin{align*}
    B := \overline{\mathrm{aco}\{f_n : n \in \mathbb{N}\}} \subseteq A.
  \end{align*} 
  Then, by $t$-orthogonality and $\inf_n \|f_n\| > 0$,
  \begin{align*}
    T : (c_0,\|\cdot\|) \to \overline{[f_1,f_2,\cdots]}, \ e_n \mapsto f_n \ (n \in \mathbb{N})
  \end{align*}
  defines a homeomorphism where $e_1, e_2, \cdots$ are the canonical unit vectors. On the other hand, since $B$ is a closed convex set, $B$ is weakly closed. Therefore, $B$ is weakly complete. Hence, $B_{c_0} = T^{-1}(B)$ is also weakly complete, which is a contradiction.
\end{proof}

The example below claims that the absolute convexity is necessary for Theorem \ref{thme1}

\begin{exam}
  Let $K$ be not locally compact, and let $X := \{e_n : n \in \mathbb\} \subseteq c_0$ where $e_1, e_2, \cdots$ are the canonical unit vectors. Then, $X$ is weakly complete, but not a norm compactoid.
\end{exam}
\begin{proof}
  Let $\lambda_1,\lambda_2, \cdots \in B_K$ such that $\inf_{i\neq j} |\lambda_i - \lambda_j| > 0$. Let us define $f \in c_0'$ by $f(e_n) = \lambda_n$ for each $n \in \mathbb{N}$. Then, we have $\inf_{i\neq j} |f(e_i) - f(e_j)| > 0$. Therefore, each Cauchy net in $A$ is eventually constant. Thus, $A$ is weakly complete, but not a norm compactoid.
\end{proof}

As in the classical case, we introduce weakly compact operators.

\begin{dfn}
  Let $(E,\|\cdot\|)$ and $(F,\|\cdot\|)$ be Banach spaces. An operator $T : E \to F$ is a weakly compact operator if $T^{**}(E'') \subseteq j_F(F)$.
\end{dfn}

\begin{prop}
  Let $(E,\|\cdot\|)$ and $(F,\|\cdot\|)$ be Banach spaces, and suppose that $(E,\|\cdot\|)$ is a norm-polar space. Let $T : E \to F$ be a continuous operator. If there exists a weakly complete subset $A$ of $F$ such that $T(B_E) \subseteq A$, then $T$ is weakly compact.
\end{prop}
\begin{proof}
  Let $\lambda \in K$ with $0 < |\lambda| < 1$. Then by Theorem \ref{ThmG1}, we have
  \begin{align*}
   \lambda B_{E''}  \subseteq \overline{j_E(B_E)}^{\sigma(E'',E')}.
  \end{align*}
  Therefore, we get
  \begin{align*}
    T^{**}(B_{E''}) &\subseteq \lambda^{-1} T^{**}\left(\overline{j_E(B_E)}^{\sigma(E'',E')}\right) \\
    &\subseteq \lambda^{-1} \overline{j_F(T(B_E))}^{\sigma(F'',F')} \\
    &\subseteq \lambda^{-1} \overline{j_F(A)}^{\sigma(F'',F')} = \lambda^{-1} j_F(A),
  \end{align*}
  which completes the proof.
\end{proof}

We recall that an operator $T \to E \to F$ is a compactoid operator if $T(B_E)$ is a norm compactoid.

\begin{thm}
   Let $(E,\|\cdot\|)$ and $(F,\|\cdot\|)$ be Banach spaces, and suppose that $(E,\|\cdot\|)$ is a norm-polar space. If $T : E \to F$ is a compactoid operator, then $T$ is a weakly compact operator.
\end{thm}
\begin{proof}
  Let $F_1 = \overline{\{0\}}^{\sigma(F,F')}$ and $(F/F_1,\|\cdot\|)$ be the quotient Banach space. Put $T_1 := \pi \circ T$ where $\pi : F \to F/F_1$ is the canonical quotient map. Then, we have $(F/F_1)'' = F''$. Therefore, it suffices to prove that $T_1$ is a weakly compact operator. \par 
  Since $T_1$ is a compactoid operator, $\overline{T_1(B_E)}$ is a complete metrizable compactoid. Thus, it follows from Proposition \ref{propjjj} that $\overline{T_1(B_E)}$ is weakly complete. Therefore, by the preceding proposition, $T_1$ is a weakly compact operator. 
\end{proof}

Of course, in the classical case, a weakly compact operator is not necessarily a compact operator. On the other hand, we have the following.

\begin{thm}
  Let $K$ be spherically complete. Then, each weakly compact operator is a compactoid operator.
\end{thm}
\begin{proof}
  Let $(E,\|\cdot\|)$ and $(F,\|\cdot\|)$ be Banach spaces, and suppose that $T : E \to F$ is a weakly compact operator. Since $B_{E''}$ is $\sigma(E'',E')$ c-compact, $T^{**}(B_{E''})$ is $\sigma(F'',F')$ c-compact. Therefore, $\overline{T^{**}(j_E(B_E))}^{\sigma(F'',F')}$ is $\sigma(F'',F')$ c-compact. Moreover, by assumption, we have $\overline{T^{**}(j_E(B_E))}^{\sigma(F'',F')} \subseteq T^{**}(B_{E''}) \subseteq j_F(F)$. Thus, we get
  \begin{align*}
    \overline{T^{**}(j_E(B_E))}^{\sigma(F'',F')} = j_F\left(\overline{T(B_E)}^{\sigma(F,F')}\right).
  \end{align*}
  Hence, we see that $\overline{T(B_E)}^{\sigma(F,F')}$ is weakly c-compact. In, particular, $\overline{T(B_E)}^{\sigma(F,F')}$ is weakly complete. Finally, by Theorem \ref{thme1}, we have that 
  \begin{align*}
    \overline{T(B_E)} = \overline{T(B_E)}^{\sigma(F,F')}
  \end{align*}
  is a norm compactoid, which completes the proof.
\end{proof}

\section{The closed range theorem} \label{clo}

If the coefficient field $K$ is spherically complete, the closed range theorem holds (\cite{close}). In this section, we prove that the closed range theorem holds between strongly norm-polar spaces (Theorem \ref{thmclo2}).

\begin{dfn}[\cite{equ}]
  A Banach space $(E,\|\cdot\|)$ is strongly norm-polar if for each closed subspace the quotient Banach space $(E/D,\|\cdot\|)$ is a norm-polar space.
\end{dfn}

\begin{rem} \label{remclo1}
  If $K$ is spherically complete and densely valued, then every Banach space is strongly norm-polar.
\end{rem}

If $(E,\|\cdot\|)$ is strongly norm-polar, then $(B_{E'},\sigma(F',F))$ is epicompact (see \cite[Section 7]{equ}). In other words, we have the following. 

\begin{thm}[{\cite[Proposition 7.1]{equ}}] \label{thmclo1}
  Let $(E,\|\cdot\|)$ be a strongly norm-polar space. Then, for each Banach space $(F,\|\cdot\|)$ and each continuous operator $T : F \to E$, we have that $(T^{*}(B_{E'}))^e$ is $\sigma(F',F)$-complete.
\end{thm}

\begin{thm} \label{thmclo2}
  Let $(E,\|\cdot\|)$ and $(F,\|\cdot\|)$ be strongly norm-polar Banach spaces, and $T : E \to F$ be a continuous operator. Then, the following are equivalent. \\
  $(1)$ $T(E)$ is norm closed. \\
  $(2)$ $T^{*}(F')$ is norm closed. \\
  $(3)$ $T^{*}(F')$ is $\sigma(E',E)$-closed. \\
  $(4)$ $\overline{T^{*}(F')} = \overline{T^{*}(F')}^{\sigma(E',E)}$.
\end{thm}
\begin{proof}
  Before beginning the proof, we note 
  \begin{itemize}
    \item strongly norm-polar Banach spaces are stable for forming of closed subspaces and quotients by closed subspaces,
    \item $\sigma \left(\overline{T(E)},(\overline{T(E)}) \right)'$ and $\sigma(F,F')$ coincide on $\overline{T(E)}$,
    \item $(E/\mathrm{Ker}T)' \subseteq E'$ is $\sigma(E',E)$-closed.
  \end{itemize} 
  Thus, we may assume that $T$ is injective and $\overline{T(E)} = F$. Then, we have that $T^{*}$ is injective. Now, the condition $(1)$ is equivalent to that $T$ is a homeomorphism. Therefore, $(1)$ implies $(3)$ and $(4)$. Clearly, $(3)$ implies $(2)$.\par
  Let us suppose $(3)$. Then, since $T$ is injective and $F'$ separates the points of $F$, we have
  \begin{align*}
    T^{*}(F') = \overline{T^{*}(F')}^{\sigma(E',E)} = E'.
  \end{align*}
  Therefore, $T^{*}$ is a homeomorphism, and so is $T^{**}$. In particular, $j_F(T(E)) = T^{**}(j_E(E))$ is closed, hence we get $(1)$. \par
  Secondly, suppose $(4)$. Then, we have 
  \begin{align*}
    \overline{T^{*}(F')} = \overline{T^{*}(F')}^{\sigma(E',E)} = E'.
  \end{align*}
    Therefore, by the same proof as that of \cite[Theorem 4]{pre}, we get $(1)$. \par
    Finally, let us suppose $(2)$. We shall prove $(3)$.  From \cite[Theorem 2.2]{kre}, it suffices to prove that $T^{*}(F') \cap B_E'$ is $\sigma(E,E')$-closed. \par
  By assumption, $T^{*}$ is a homeomorphism onto its image. Therefore, there exists $\lambda \in K$, $\lambda \neq 0$, such that 
    \begin{align*}
      \|T^{*}(f)\| \ge |\lambda| \|f\| \ \mathrm{for} \ \mathrm{each} \ f \in F'.
    \end{align*}
  Thus, we get the equality
  \begin{align*}
    T^{*}(F') \cap B_E' = T^{*}(B_F'(0,1/\lambda)) \cap B_{E'}.
  \end{align*}
  Therefore, by Theorem \ref{thmclo1},
  \begin{align*}
    T^{*}(F') \cap B_E' &= (T^{*}(F') \cap B_E')^e \\
    &= (T^{*}(B_F'(0,1/\lambda)) \cap B_{E'})^e \\
    &= (T^{*}(B_F'(0,1/\lambda))^e \cap B_{E'}
  \end{align*}
  is $\sigma(E',E)$-closed, which completes the proof.
\end{proof}

\end{document}